\documentclass[12pt]{article}
\usepackage[utf8]{inputenc}
\usepackage[T1]{fontenc}
\usepackage{lmodern}
\usepackage{graphicx}
\newcommand{\vs}{\vskip10pt}

\usepackage{amsmath}
\usepackage{amsthm}
\usepackage{amssymb}
\usepackage[export]{adjustbox}
\usepackage[margin=0.75in]{geometry}
\usepackage{setspace}
\usepackage{enumitem}
\usepackage{float}
\usepackage{tabularx,ragged2e}
\newtheorem{thm}{Theorem}[section]
\newtheorem{defn}[thm]{Definition}
\newtheorem{prop}[thm]{Proposition}

\newtheorem{lem}[thm]{Lemma}
\newtheorem*{lem*}{Lemma}
\newtheorem*{thm*}{Theorem}
\newtheorem*{cor*}{Corollary}
\newtheorem*{rem*}{Remark}

\newtheorem*{clm*}{Claim}

\newtheorem*{ques*}{Question}
\usepackage[all]{xy}
\usepackage{tikz-cd}
\usepackage{mathtools}
\usepackage{comment}
\usepackage{authblk}
\usepackage{url}
\makeatletter
\usepackage{etoolbox}

\AtEndDocument{%
	\par
	\newpage
	\medskip
	\begin{tabular}{@{}l@{}}%
		\textsc{Chris Karpinski}\\
		Department of Mathematics and Statistics, McGill University, Burnside Hall, \\
		805 Sherbrooke Street West, \\
		Montreal, QC, H3A 0B9, Canada. \\
		\textit{E-mail address}: \texttt{christopher.karpinski@mcgill.ca}
\end{tabular}}
\AtEndDocument{%
	\par
	\medskip
	\begin{tabular}{@{}l@{}}%
		\textsc{Damian Osajda}\\
		Department of Mathematical Sciences,
		University of Copenhagen,
		\\Nørregade 10, 1165 København, Denmark.
        \\
        Institute of Mathematics, Polish Academy of Sciences,
	    \\\'Sniadeckich 8, 00-656 War\-sza\-wa, Poland. \\
		\textit{E-mail address}: \texttt{dosaj@math.uni.wroc.pl}
\end{tabular}}
\AtEndDocument{%
	\par
	\medskip
	\begin{tabular}{@{}l@{}}%
		\textsc{Piotr Przytycki}\\
		Department of Mathematics and Statistics, McGill University, Burnside Hall, \\
		805 Sherbrooke Street West, \\
		Montreal, QC, H3A 0B9, Canada. \\
		\textit{E-mail address}: \texttt{piotr.przytycki@mcgill.ca}
\end{tabular}}

\title{Weak Tits alternative for \textcolor{black}{groups acting geometrically on} buildings}
\date{}
\author{Chris Karpinski, Damian Osajda, and Piotr Przytycki}

\usepackage{lipsum}
\newcommand{\R}{\mathbb{R}}
\newcommand{\N}{\mathbb{N}}

\newcommand{\Z}{\mathbb{Z}}
\newcommand{\pcom}[1]{\textcolor{black}{#1}}

\usepackage{hyperref}
\hypersetup{
	colorlinks=true,
	linkcolor=red,
	citecolor=blue,
}

\begin{document}
	
	\maketitle

	\begin{abstract}
		
		We show that if a group $G$ acts geometrically by type-preserving automorphisms on a building, then $G$ satisfies the weak Tits alternative, namely, that $G$ is either virtually abelian or contains a non-abelian free group. 
		
	\end{abstract}
	
	\section{Introduction}
	
	Buildings were introduced by Jacques Tits in the 1950s as a tool to study semisimple algebraic groups. Since their inception, buildings have found diverse applications throughout mathematics, well beyond their roots in the theory of algebraic groups; see for instance the survey article \cite{BuildingsApp}. 
	
	An ongoing area of interest has been the study of algebraic properties of groups acting on buildings.
    Buildings might be equipped with a structure of a nonpositively curved metric space (see e.g.\ \cite[Chapter 18]{Davis}), so it is believed that groups acting on them in a nice enough manner exhibit a property shared by many ``non-positively curved'' groups: the \emph{Tits alternative}. The Tits alternative is a dichotomy for groups and their subgroups, \textcolor{black}{first formally defined by Bass--Lubotzky \cite[Remark~1.4(2)]{MR0693651} but whose initial study dates back to the work of Jacques Tits in \cite{Tits1}}, where it was shown that every finitely generated linear group is either virtually solvable or contains the free group $F_2$ as a subgroup. We will consider a weaker version of the Tits alternative: we will say that a group satisfies the \emph{weak Tits alternative} if it is either virtually abelian or contains $F_2$ as a subgroup. The weak Tits alternative has been shown to be satisfied for groups acting properly and cocompactly on Euclidean buildings in \pcom{\cite[Theorem~F]{BallmanBrin} (where in the virtually abelian case the group is finite-by-Bieberbach)}. 
    Sageev and Wise show in \cite{SageevWise} that groups \textcolor{black}{with a bound on the cardinality of finite subgroups acting properly on finite-dimensional CAT(0) cube complexes satisfy the Tits alternative}. In particular, this implies the Tits alternative for such groups acting properly on right-angled buildings.
    The Tits alternative was proved for groups \textcolor{black}{with a bound on the cardinality of finite subgroups acting properly on 2-dimensional
    complexes with some ``non-positive curvature'' features in \cite{MR4305241,osajda_przytycki_2022}}. This covers the case of all 2-dimensional buildings.
    
	In this paper, we extend the above results obtained for Euclidean, right-angled, and 2-dimensional buildings to actions on arbitrary buildings. Our main theorem is the following:

    \newtheorem{theorem}{Theorem}
	\renewcommand*{\thetheorem}{\Alph{theorem}}
	\newtheorem{corollary}[theorem]{Corollary}
	\renewcommand*{\thecorollary}{\Alph{theorem}}
	
	\begin{theorem}
		\label{A}
		Let $G$ be a group acting properly and cocompactly (i.e.\ geometrically) by type-preserving automorphisms on a building. Then $G$ is either virtually abelian or contains a non-abelian free subgroup. 
		
	\end{theorem}

\pcom{Theorem~\ref{A} still holds if instead of requiring that each automorphism $g\in G$ is type-preserving, we assume that $g$ descends to a permutation of the set of types, which follows from Theorem~\ref{A} applied to a finite index subgroup. However, we do not know if Theorem~\ref{A} still holds if $G$ contains an automorphism that does not descend to such a permutation.}

\medskip 

	\noindent \textbf{Proof outline.} Our proof consists of first removing a possible finite factor of the underlying Coxeter group $W$ of the building (Lemma \ref{3.9}) and then splitting into the cases of whether or not the building is thin. In the case of the building being thin, the weak Tits alternative for $G$ follows quickly by purely algebraic arguments from the classical Tits alternative for linear groups.

	In the non-thin case, our proof relies on the construction of a tree of chambers in the building and group elements $g,g' \in G$ acting on this tree. Our construction relies on probabilistic arguments adapted from the proof of \cite[Lemma 2.10]{MR4305241}, originally stemming from arguments in \cite{BallmanBrin}. More precisely, \textcolor{black}{the} basic idea of the proof in the non-thin case is the following.

	We begin with a branching panel in some wall $\Omega$ in an apartment $\Sigma$ of a building $\Delta$. By Lemma~\ref{3.8}, there exists a wall $\Omega'$ in $\Sigma$ which is parallel to $\Omega$. We then connect $\Omega$ to $\Omega'$ via a minimum length gallery $\gamma$ between pairs of panels on these walls. Let $\sigma$ be the  panel in $\Omega$ containing the initial chamber of $\gamma$. By Lemma \ref{3.5}, we have that every panel in $\Omega$ branches, so that $\sigma$ branches. 
 
    Using Proposition \ref{3.6} applied to pairs of three chambers in $\sigma$, we produce a \color{black} ``tree of chambers'' and group elements $g,g' \in G$ which act on this tree. By examining the orbit of $\sigma$ under $\langle g,g' \rangle$, \pcom{we show that $g$ and $g'$ generate a free subgroup of $G$}.  \color{black}
	
	\vs 

    \color{black}
    
    Pierre-Emmanuel Caprace has pointed out to us that Theorem \ref{A} follows from the following more general result on the weak Tits alternative for groups acting geometrically on products of $\mathrm{CAT(0)}$ spaces satisfying rank rigidity.

\begin{theorem}
    \label{B}
    Let $G$ be a group acting geometrically on a product $X=X_1 \times \cdots \times X_n$ of $\mathrm{CAT(0)}$ spaces $X_1,\ldots,X_n$ such that no $X_i$ is a point and each $X_i$ satisfies rank rigidity. In the case that each $X_i$ is a Euclidean building, assume that the action is type-preserving. Then $G$ satisfies the weak Tits alternative. 
\end{theorem}

In the appendix, we \pcom{outline} a proof of Theorem \ref{B} (due to Caprace) and we show how Theorem~\ref{B} implies Theorem \ref{A}.

\vs



 \color{black}

	\vs

\noindent\textbf{Acknowledgements}: We thank Pierre-Emmanuel Caprace for pointing out to us an alternate proof of our main theorem using rank rigidity (outlined in the appendix). We also thank Chris Hruska for bringing a result of Wilking to our attention that simplified the proof of Proposition~\ref{3.10}, and we thank Marcin Sabok for suggesting the proof of Lemma~\ref{3.51}. \textcolor{black}{We also thank the referees for a thorough reading of the article and for many insightful corrections and suggestions.} DO and PP were partially supported by (Polish) Narodowe Centrum Nauki, UMO-2018/30/M/ST1/00668. PP was partially supported by NSERC.
	
	\section{Preliminaries}
	
	\subsection{Chamber systems} 
	
	The following definitions are from \cite{Ronan}. 
	
	\begin{defn}
		
		A \textbf{chamber system} is a set $C$ together with a set $I$ such that each element~$i$ of~$I$ determines a partition of $C$. Two elements in the same part of $C$ determined by $i \in I$ are called \textbf{$i$-adjacent}, and we will call two elements of $C$ \textbf{adjacent} if they are $i$-adjacent for some $i \in I$. The elements of $C$ are called \textbf{chambers} and we refer to $I$ as the \textbf{index set}. 
		
	\end{defn}
	
	A \textbf{gallery} is a finite sequence of chambers $(c_0,\ldots,c_k)$ such that each $c_{j-1}$ is adjacent to $c_j$ and $c_{j-1} \neq c_j$. A \textbf{subgallery} of a gallery $(c_0,\ldots, c_k)$ is a subsequence of $(c_0, \ldots, c_k)$ of the form $(c_i, c_{i+1},\ldots,c_j)$ for some $0 \leq i \leq j \leq k$. Given a gallery $\gamma = (c_0,\ldots, c_k)$, the \textbf{inverse gallery} is the gallery $\gamma^{-1} := (c_k, c_{k-1}, \ldots, c_0)$. The gallery $(c_0,\ldots,c_k)$ has \textbf{type} $i_1\cdots i_k \in I^*$ (where $I^*$ denotes the set of all finite length words in elements of $I$) if $c_{j-1}$ is $i_j$-adjacent to $c_j$. The \textbf{length} of a gallery~$\gamma$, denoted $\ell(\gamma)$, is the length of its type as a word in $I^*$. A \textbf{geodesic gallery} is a gallery that has minimal length among all galleries with the same initial and terminal chambers. If each $i_j$ belongs to a fixed subset $J \subseteq I$, then we call the gallery $(c_0,\ldots,c_k)$  a \textbf{$J$-gallery}.

	A chamber system $C$ over a set $I$ is called \textbf{connected} (resp.\ \textbf{$J$-connected}) if any pair of chambers can be joined by a gallery (resp.\ $J$-gallery). The $J$-connected components are called \textbf{$J$-residues}. For $i \in I$, an $\{i\}$-residue is called a \textbf{panel}, whose \textbf{type} is $i$.  If $\sigma$ is a panel, we will say that each chamber $c \in \sigma$ \textbf{has} $\sigma$ as a panel. By a gallery between panels $\alpha, \sigma$, we mean a gallery between a chamber  in $\alpha$ and a chamber in $\sigma$.  
	
	
	\subsection{Coxeter groups}
	
	\begin{defn}
		
		A \textbf{Coxeter group} is a group $W$ having a \textbf{Coxeter presentation}, that is, a presentation of the form: 
		
		$$ W = \langle S \vert s^2 = 1 = (rs)^{m_{rs}} \text{ for all } r \neq s \text{ in } S , m_{rs} \in \{2,3,\ldots ,\infty\} \text{ and } m_{rs} = m_{sr} \rangle $$
		
		where $m_{rs} = \infty$ means that there is no relation between $r,s$. 
	\end{defn}

	Given a Coxeter presentation as above, we say that $S$ is a \textbf{Coxeter generating set} of $W$ and that $(W,S)$ is a \textbf{Coxeter system}. 
 A conjugate $s^w := wsw^{-1}$ for $w \in W$ and $s \in S$  is called a \textbf{reflection}.

	Given a Coxeter system $(W,S)$, we denote by $\vert \cdot \vert_S$ the word length of an element of $W$ with respect to $S$ (i.e.\ for $w \in W$, $\vert w \vert_S$ represents the length of \textcolor{black}{a} shortest word over $S$ representing $w$ in $W$) and we denote $d_S$ the word metric on $W$ with respect to $S$ (i.e.\ $d_S(u,v) = \vert u^{-1}v \vert_S$ for each $u,v \in W$). A \textbf{standard Coxeter subgroup} of a Coxeter group $W$ with Coxeter generating set $S$ is a subgroup generated by some $T \subseteq S$. \textcolor{black}{This terminology is justified by the fact that $(\langle T \rangle, T)$ is a Coxeter system for each subset $T \subseteq S$.}

	A key example of a chamber system is a Coxeter system $(W,S)$. Here, the set of chambers is~$W$ and the index set is  $S$. Two chambers $w_1, w_2 \in W$ are $s$-adjacent for $s \in S$ if $w_2 = w_1 s$ in $W$. 
	
	\subsection{Buildings}
	
	For background on buildings, we follow the books \cite{BrownAbr} and \cite{Ronan}. 
	
	\begin{defn}
		
		A \textbf{building} of type $(W,S)$ is a chamber system $\Delta$ over $S$ such that each panel contains at least two chambers, equipped with a map $\delta : \Delta \times \Delta \to W$ such that if $f$ is a geodesic word over $S$, then $\delta(x,y) = f \in W$ if and only if $x,y$ can be joined by a gallery of type $f$. The map~$\delta$ is called a \textbf{$W$-metric on $\Delta$}. 
		
	\end{defn}

	Note that in a building $\Delta$, two adjacent chambers $x,y$ are $s$-adjacent for a unique $s \in S$ ($s = \delta(x,y)$). We also have a metric $d$ on $\Delta$ defined by $d(x,y) = \vert \delta(x,y)\vert_S$ for each $x,y \in \Delta$. We will refer to $d$ as the \textbf{gallery metric} (note that the triangle inequality for $d$ follows from \cite[Lemma 5.28]{BrownAbr}, so $d$ is indeed a metric). We will use the notation $(\Delta, \delta)$ to denote a building $\Delta$ with its associated $W$-metric $\delta$. 
	
	A \textbf{type-preserving automorphism} $\phi$ of a building $(\Delta,\delta)$ is a bijective map $\phi: \Delta \to \Delta$ that preserves the $W$-metric $\delta$, i.e.\ $\delta(x,y) = \delta(\phi(x), \phi(y))$ for each $x,y \in \Delta$.  
	
	Given a panel $\sigma$ in a building, the \textbf{degree} of $\sigma$, denoted $\deg \sigma$, is the number of chambers in the building having $\sigma$ as a panel.  A panel is \textbf{branching} if it has degree at least 3. A building is called \textbf{thin} if it has no branching panels, i.e.\ each panel has degree 2 and hence is a panel of exactly two chambers.

	Given a building $\Delta$ of type $(W,S)$ with $W$-metric $\delta$, a subset $\Delta_2 \subseteq \Delta$ is  a \textbf{subbuilding} if $(\Delta_2, \delta\vert_{\Delta_2})$ is a building (of possibly different type than $\Delta$). 

	Note that a Coxeter group $W$ is an example of a building, where we take $W$ as the set of chambers and a Coxeter generating set $S$ as the index set, and equip $W$ with $W$-metric $\delta_W$ defined by $\delta_W(x,y) = x^{-1}y$ for each $x,y \in W$. The corresponding gallery metric is the word metric $d_S$ on~$W$.

	A Coxeter group $W$ admits a type-preserving action on its associated building induced from the  action on itself by left multiplication. When viewing a Coxeter group $W$ as a building, we have the notion of \textbf{walls} that separate $W$  into two connected components. Given a reflection $r=s^u \in W$ (for $s \in S$ and $u \in W$), the \textbf{wall} associated to $r$ is the set  $M_r = \{\{c_1, c_2\} : c_1, c_2 \in W \text{ are adjacent and } rc_1 = c_2\}$.  Thus, a wall is the set of all panels fixed by $r$. The sets $\alpha_r^+ = \{w \in W : d_S(w,u) < d_S(w, us)\}$ and  $\alpha_r^- = \{w \in W : d_S(w,u) >  d_S(w, us)\}$ are called the \textbf{roots} of $r$. Two walls $M_r$ and $M_s$ are \textbf{parallel} if $\langle r,s \rangle \cong D_{\infty}$. 
	
	Buildings contain special subspaces, called \textbf{apartments}. Let $(W,S)$ be a Coxeter system and let $\Delta$ be a building of type $(W,S)$. For a subset $X \subseteq W$, a map $\alpha : X \to \Delta$ is an \textbf{isometry} if it preserves $W$-distance: $\delta_{\Delta}(\alpha(x), \alpha(y)) = \delta_W(x,y)$ for each $x,y \in X$. An \textbf{apartment} in a building~$\Delta$ of type $(W,S)$ is an image $\alpha(W)$ for an isometry $\alpha : W \to \Delta$.
	
	By the characterization in \cite[Section 5.5.2]{BrownAbr}, the apartments of a building $\Delta$ are precisely the thin subbuildings of $\Delta$.  Every panel in an apartment is contained in a panel in the building. For a panel $\sigma$ in an apartment, we say that $\sigma$ is \textbf{branching} if the panel in the building containing $\sigma$ is branching. 
	
	A \textbf{wall} (resp.\ \textbf{root}) in an apartment of $\Delta$ is an isometric image of a wall (resp.\ root) in $W$. A gallery $\alpha$ \textbf{crosses} a wall $\Omega$ if $\alpha$ passes through both chambers in a panel of $\Omega$. 
	
	\section{Proof of the main theorem}
	
	Recall that an action of a group $G$ on a metric space $(X,d)$ is called $\textbf{geometric}$ if it is \textbf{proper} (i.e.\ for each $x \in X$, $r \geq 0$, we have $\vert \{g \in G : d(x,gx) \leq r\} \vert < \infty$) and \textbf{cocompact} (i.e.\ there is a compact fundamental domain for the action $G \curvearrowright X$). We fix a group $G$ acting geometrically by type-preserving automorphisms on $(\Delta,d)$, where $(\Delta, \delta)$ is a building of type $(W,S)$, with $d$ the gallery metric. 
 \pcom{Note that since $G \curvearrowright (\Delta,d)$ is cocompact, 
there are only finitely many $G$-orbits of chambers. By the properness of the action of $G$, each ball of finite radius contains only a finite subset of each $G$-orbit. Consequently, the metric space $(\Delta, d)$ is locally finite (i.e.\ balls of finite radius are finite). In particular, $|S|$ is finite.}

	We consider two cases on the building $\Delta$: the case of $\Delta$ being thin (equivalently, consisting of a single apartment) and the complementary case of $\Delta$ not being thin, hence consisting of more than one apartment. 
	
	\begin{prop}
		
		\label{3.10}
		
		Let $(\Delta, \delta)$ be a thin building and let $G$ act geometrically on $\Delta$ by type-preserving automorphisms. Then $G$ is either virtually abelian or contains $F_2$ as a subgroup.
		
	\end{prop}

	\begin{proof}
		
		Since $\Delta$ is thin, it consists of a single apartment.  Since $G$ acts by type-preserving automorphisms of $\Delta$ and since $W$ is isomorphic to the group of all type-preserving automorphisms of $\Delta$ (by \cite[Proposition 3.32]{BrownAbr}), we have a group homomorphism $\rho : G \to W$, given by fixing a chamber $c$ and putting $g \mapsto \delta(c, gc)$. Since the action of $G$ on $\Delta$ is proper, we have that the stabilizer $\mathrm{Stab}(c)$ is finite and hence that $\ker \rho$ is finite. Also, since the action of $G$ on $\Delta$ is cocompact, we have that $\rho(G)$ is of finite index in $W$. 
		
		Since $W$ is linear over $\R$ (see, for instance, \cite[Section 2.5]{BrownAbr}), we have that $W$  satisfies the classical Tits alternative: every subgroup of $W$ is either virtually solvable or contains $F_2$. Therefore, $\rho(G)$ is either virtually solvable or contains $F_2$. \textcolor{black}{It is also proved in \cite{Noskov2002} that $W$ satisfies a strong version of the classical Tits alternative.}
		
		If $\rho(G)$ contains a subgroup $H \cong F_2$, then $\rho^{-1}(H) \leq G$ surjects onto $H \cong F_2$, and hence contains $F_2$, so $G$ contains $F_2$.
		
		If $\rho(G)$ is virtually solvable, then $\rho(G)$ is virtually abelian, since Coxeter groups are CAT(0) (since they act geometrically on their Davis complex, which is a CAT(0) space; see Chapters 7 and 12 of \cite{Davis}), and solvable subgroups of CAT(0) groups are virtually abelian by \cite[Theorem III.$\Gamma$.1.1(3)]{BH99}. Therefore, letting $H < \rho(G)$ be a finite-index free abelian subgroup, we have that $\rho^{-1}(H)$ has finite index in $G$ and surjects onto a finitely generated free abelian group with finite kernel. By \cite[Theorem 2.1 $c) \implies d)$]{MR1783960}, we have that $\rho^{-1}(H)$ is virtually abelian, hence $G$ is virtually abelian. 
		
		\end{proof}

	We now move on to the case where $\Delta$ is not thin. 
	
	\begin{prop}
		
		\label{3.11}
		
		Let $(\Delta, \delta)$ be a building of type $(W,S)$ that is not thin and such that $W$ does not decompose as $W \cong W_1 \times W_2$ for $W_1, W_2$ standard Coxeter subgroups of $W$, with $W_1$ finite and non-trivial. Let $G$ act geometrically on $\Delta$ by type-preserving automorphisms. Then $G$ contains $F_2$ as a subgroup. 
		
	\end{prop}

	In the proof of  Proposition \ref{3.11}, we will need the following lemma. It was first stated in \cite[Lemma 4.1]{DavisShapiro}, and later in \cite{BipolarCox}, where a different proof was given. The proof in \cite{BipolarCox} relies on \cite[Lemma 8.2]{TwistRig} and the parallel wall theorem (\cite[Theorem 2.8]{BrinkB}).

	\begin{lem}
		
		\label{3.8}
		
		If $(W,S)$ is a Coxeter system such that $W$ does not decompose as the direct product of standard Coxeter subgroups $W_1$, $W_2$, where $W_1$ is a finite non-trivial Coxeter group, then for each wall $\Omega$ in $W$, there exists a wall $\Omega'$ in $W$ which is parallel to $\Omega$. 
		
	\end{lem}

	The following lemma and Lemma \ref{3.8} allow us to reduce to the case where we can find a wall disjoint from any given wall. 

	\begin{lem}
		
		\label {3.9}
		
		If $G$ acts geometrically by type-preserving automorphisms on a building $(\Delta, \delta)$ of type $(W,S)$ and $W \cong W_1 \times W_2$, where $W_1 = \langle S_1 \rangle$ and $W_2 = \langle S_2 \rangle $ are standard Coxeter subgroups of $W$ with $S_1 \coprod S_2 = S$ and with $W_1$ finite,  then there exists a building $\Delta_2$ of type $(W_2,S_2)$ on which $G$ acts geometrically by type-preserving automorphisms. 
		
	\end{lem}

	\begin{proof}
		

        \color{black}
        
        By \cite[Theorem 3.10]{Ronan}, we have a decomposition $\Delta \cong \Delta_1 \times \Delta_2$, where each $\Delta_i$ is a building of type $(W_i, S_i)$. Since the action of $G$ on $\Delta$ is type-preserving, this product decomposition is preserved by $G$. Hence, the projection map $\Delta \to \Delta_2$ is $G$-equivariant, and so $G$ acts on $\Delta_2$ cocompactly by type-preserving automorphisms. Since $\Delta_1$ is locally finite and $W_1$ is finite, it follows that $\Delta_1$ is finite and so the properness of $G \curvearrowright \Delta$ implies that $G$ acts on $\Delta_2$ properly. Therefore, $G$ acts on~$\Delta_2$ geometrically by type-preserving automorphisms. 

        \color{black}

		\end{proof}

	Combining the results of Proposition \ref{3.10}, Proposition \ref{3.11} and Lemma \ref{3.9}, we conclude the proof of Theorem \ref{A}. It remains to prove Proposition \ref{3.11}. 


	\section{Proof of Proposition \ref{3.11}}
	
	In the proof of Proposition \ref{3.11}, we will need the following lemma. 
	
	\begin{lem}[\textcolor{black}{\cite[Corollary 4.2]{MR2110465}}]
		
		\label{3.5}
		
		Let $\Omega$ be a wall in an apartment $\Sigma$ of a building $(\Delta, \delta)$. Suppose that $\Omega$ has a branching panel. Then every panel of $\Omega$ branches. 
		
	\end{lem}

	For the remainder of this section, we fix an apartment $\Sigma$ containing a branching panel. Let $\Omega \subseteq \Sigma$ be any wall containing this branching panel. Invoking Lemma \ref{3.8}, there exists a wall $\Omega'$ \pcom{in $\Sigma$} which is parallel to $\Omega$. 
 Fix a geodesic gallery $\gamma$ which has minimal length among all galleries joining a chamber inside a panel of $\Omega$ and a chamber inside a panel of $\Omega'$. Let $s_0$ be the type of the panel in $\Omega$ containing the first chamber of $\gamma$ and let $s_k$ be the type of the panel in $\Omega'$ containing the last chamber of $\gamma$. Let $s_1\cdots s_{k-1}$ be the type of $\gamma$, \pcom{where} $s_i \in S$. For $j = 1,\ldots ,k-1$, put $s_{k+j} = s_{k-j}$.

	\begin{lem}
		
		\label{3.2}
		
		Let $w = s_0\cdots s_{2k-1} \in W$. Then for any $n \in \N$, we have that \newline$s_0w^n = s_1s_2 \cdots s_{2k-1} (s_0 s_1 \cdots s_{2k-1})^{n-1}$ is a geodesic word in $W$ (by convention, $s_0$ cancels the first letter $s_0$ of $w^n$).  
		
	\end{lem}
	
	\begin{proof}
		
		Let $s,r \in W$ be such that $\Omega = M_r, \Omega' = M_s$.  Let $\gamma$ be the above  minimal length geodesic gallery between panels in the walls $\Omega$ and $\Omega'$  having type $s_0 w = s_1\cdots s_{k-1}$.  Then $s_0w^n = s_1s_2 \cdots s_{2k-1} (s_0s_1 \cdots s_{2k-1})^{n-1}$ is the type of the gallery $\gamma_n := \bigcup_{i=0}^{n-1} (sr)^i (\gamma \cup s \gamma)$ (see Figure~\ref{fig:geodesics-concat2}).

\begin{figure}[H]
	\centering
	\includegraphics[width=0.8\linewidth]{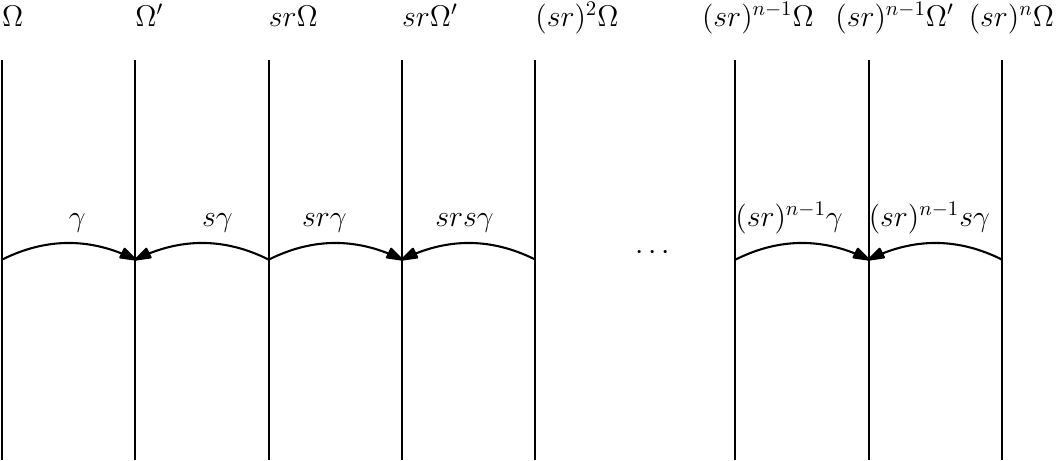}
	\caption{The concatenation of \pcom{the geodesic galleries forming $\gamma_n$.}}
	\label{fig:geodesics-concat2}
\end{figure}

		We claim that $\gamma_n$ is a geodesic gallery. Indeed, if $\alpha$ were a gallery with the same starting and ending chambers as $\gamma_n$, then by \cite[Lemma 2.5(ii)]{Ronan}, we must have that $\alpha$ crosses each such wall $(sr)^i \Omega$. Let $\alpha_i$ denote the segment of $\alpha$ between the successive parallel walls $(sr)^i\Omega$ and $(sr)^{i+1}\Omega$. Then $(sr)^{-i} \alpha_i$ is a gallery between the walls $\Omega$ and $(sr)\Omega = s \Omega$. Since the chambers contained in $\Omega$ and $s \Omega$ are in different roots of $s$, by \cite[Lemma 2.5(ii)]{Ronan} we have that $(sr)^{-i} \alpha_i$ crosses $\Omega'$.  Since $\gamma$ is a minimum length geodesic gallery between panels in the walls $\Omega$ and $\Omega'$, we have that $s \gamma^{-1}$ is a minimum length gallery between panels in $\Omega'$ and $s \Omega$. Therefore, the length of the subgallery of $(sr)^{-i} \alpha_i$ between $\Omega$ and $\Omega'$ is at least $\ell(\gamma)$ and similarly the  length of the subgallery of $(sr)^{-i} \alpha_i$ between $\Omega'$ and $s \Omega$ is at least $\ell(s\gamma^{-1})$. Therefore, $\ell((sr)^{-i} \alpha_i) \geq \ell(\gamma \cup s \gamma^{-1})$. Translating by $(sr)^i$, we obtain $\ell(\alpha_i) \geq \ell((sr)^i (\gamma \cup s \gamma^{-1}))$. As this holds for all $i = 0,1,\ldots, n-1$, we conclude that $\ell(\alpha) \geq \ell(\gamma_n)$, and therefore that $\gamma_n$ is a geodesic gallery. 
		Since $s_0w^n$ is the type of $\gamma_n$, it follows that $s_0w^n$ is a geodesic word. 
		\end{proof}
	
	Using ideas from the work of Ballmann and Brin in \cite{BallmanBrin}, we construct the following Markov chain. The set $A$ of states will consist of pairs $(c, j)$, where $c$ is a chamber of $\Delta$ and $j \in \Z$ is an index taken modulo $2k$ (recall that $k$ is the length of \textcolor{black}{a} minimum length gallery $\gamma$ between the walls $\Omega$ and $\Omega'$). \textcolor{black}{The role of the indices $j$ is to identify a suitable notion of adjacency of states, which we use below to define the transition probability in our Markov chain. As we will see below, we want states to be adjacent in the Markov chain if they consist of chambers sharing a panel of the ``correct type'', determined by their indices. This will ultimately allow us to define the appropriate branches of the tree of chambers that we will construct inside $\Delta$ to produce our desired $F_2$ subgroup in $G$. }
 
    We define the \textbf{transition probability} $p(a \to a')$ from $a = (c,j) \in A$ to $a' = (c',i) \in A$ to be positive if $i = j +1$ and $c$ and $c'$ share a panel $\sigma$ with type $s_j$,  in which case we set $p(a \to a') =  \frac{1}{\deg \sigma - 1}$, otherwise we put $p(a \to a') = 0$. \textcolor{black}{Note that $\deg \sigma$ is finite since $\Delta$ is a locally finite building}. We have an action of $G$ on $A$ via $g  (c,j) = (g  c, j)$ for each chamber $c$ of $\Delta$ and $j = 0,1,\ldots ,2k-1$. 
	
	Given a sequence of states $a_0,\ldots,a_n$, we put $p_n(a_0,\ldots,a_n) = \prod_{i=0}^{n-1} p(a_i \to a_{i+1})$. Given a finite sequence $(a_0,\ldots,a_n)$ of states, \pcom{and $N\in \Z$}, denote the \textbf{cylinder set} $[a_0\cdots a_n]_{(N,N + n)} := \{(b_i)_{i \in \Z} : b_{i} = a_{i-N} \text{ for all } i = N,\ldots,N+n\} \subseteq A^{\Z}$.  
	
	\begin{lem}
		
		There exists a \textcolor{black}{unique} shift-invariant measure $\mu$ on $A^{\Z}$ such that $\mu([a_0\cdots a_n]_{(N, N+n)}) = p_{n}(a_0,\ldots ,a_n)$ for each $a_0,\ldots ,a_n \in A$ ($n \geq 0$), \pcom{and $N\in \Z$}. 
		
	\end{lem}

	\begin{proof}
		
		By \cite[Example (8)]{Walters}, we need to check that the following properties of $p$ are satisfied: 
		
		\begin{enumerate}[label = (\roman*)]
			\item For any $a \in A$, $\sum_{a' \in A} p(a \to a') = 1$
			\item For any $a \in A$, $\sum_{a' \in A} p(a' \to a) = 1$
		\end{enumerate}
		
		For (i),  given $a = (c,j) \in A$, we have $p(a \to a') \neq 0$ only if the chambers of $a'$ and $a$ share a panel $\sigma$ of type $s_{j}$. In this case, we then have $p(a \to a') = \frac{1}{\deg \sigma - 1}$. Since there are exactly $\deg \sigma - 1$ chambers other than the chamber of $a$ having $\sigma$ as a panel, we obtain: 
		
		$$\sum_{a' \in A} p(a \to a') = (\deg \sigma - 1) \cdot \frac{1}{\deg \sigma - 1} = 1$$
		
		For (ii), given $a = (c,j) \in A$, we have $p(a' \to a) \neq 0$ only if the chambers of $a'$ and $a$ share a panel $\sigma$ of type $s_{j-1}$, and in this case we have $p(a' \to a) = \frac{1}{\deg \sigma -1}$. We then have: 
		
		$$ \sum_{a' \in A} p(a' \to a) = (\deg \sigma - 1) \cdot \frac{1}{\deg \sigma - 1} = 1$$
		
		Therefore, $p$ induces a shift invariant measure $\mu$ on $A^{\Z}$ with the desired value on cylinder sets. 
		\end{proof}

	The measure $\mu$ is $G$-invariant, since the action of $G$ on $\Delta$ is type-preserving and hence preserves adjacency. \textcolor{black}{The measure $\mu$ can be thought of as an analogue of the Liouville measure on the geodesic flow in the case of manifolds,} \pcom{where certain flow directions are distinguished.}
 \color{black}
 
 \pcom{Below, the Borel structure on $A^{\Z}$ comes from putting the discrete topology on $A$ and equipping~$A^{\Z}$ with the product topology.}
 
 \begin{lem}
     \label{3.51}
     There is a measurable fundamental domain $D$ for the action of $G$ on $A^{\Z}$. In particular, the measure $\mu$ on $A^{\Z}$ descends to a measure $\bar{\mu}$ on $A^{\Z}/G$ defined as the restriction of $\mu$ to $D$. Moreover, $\bar{\mu}$ is shift invariant and finite. 
 \end{lem}

 \begin{proof}
     By \cite[Theorem 6.4.4]{Gao} and the fact that $G$ is countable (since $G \curvearrowright (\Delta,d)$ is proper and $\Delta$ is countable), it suffices to show that the orbit equivalence relation $E_G$ of $G \curvearrowright A^{\Z}$ is a closed subset of $A^{\Z} \times A^{\Z}$. Let $((a^n)_n, (g_na^n)_n)$ be a sequence in $E_G$ converging to $(a^{\infty}, b^{\infty}) \in A^{\Z} \times A^{\Z}$. Then the sequence $(a_0^n)_n$ is eventually equal to $a_0^{\infty}$ and the sequence $(g_na_0^n)_n$ is eventually equal to~$b_0^{\infty}$. Thus, we have that $g_n a_0^{\infty} = b_0^{\infty}$ for all sufficiently large $n$. By \pcom{the} properness of the action of $G$ on $\Delta$, there are only finitely many $g \in G$ with $ga_0^{\infty} = b_0^{\infty}$, so passing to a subsequence, we may assume that the sequence $(g_n)_n$ is constant, with all $g_n$ equal to some $g \in G$. We then have that $b^{\infty} = \lim_{n \to \infty}(ga^n)_n= ga^{\infty}$, so we conclude that $(a^{\infty}, b^{\infty}) \in E_G$ and hence that $E_G$ is closed.
	
	Since $\mu$ is invariant under the forward shift map $T : A^{\Z} \to A^{\Z}$ given by $T((a_i)_{i \in \Z}) = (a_{i+1})_{i \in \Z}$, so is $\bar{\mu}$. 
	
	Since the action of $G$ on $\Delta$ is cocompact, we have that $A / G$ is finite, and so the measure $\bar{\mu}$ is finite. 
 \end{proof}
 \color{black}
	To produce an $F_2$ subgroup inside $G$, we will construct a tree of chambers in $\Delta$ and a subgroup of $G$ acting freely on this tree using the Poincaré recurrence lemma. The next proposition is the key ingredient involved in this construction. 
	
	\begin{prop}
		
		\label{3.6}
		
		Given a pair of states $(d, 0)$ and $(d',1)$ with $p((d,0) \to  (d',1)) > 0$, there exists a sequence $(a_1,\ldots ,a_n)$ of elements of $A$ with the following properties:
		
		\begin{itemize}
			\item $a_1 = (d',1)$  and $a_n = g(d,0)$ for some $g \in G$.
			\item $p(a_i \to a_{i+1}) > 0$ for each $i = 0,1,\ldots ,n-1$.
		\end{itemize}
	
	\end{prop}
	
	\begin{proof}
		
		Denoting $a = (d,0)$ and $a' = (d',1)$, consider the cylinder set $[aa'] _{(0,1)} = \{(a_i)_{i \in \Z} : a_0 = a, a_1 = a'\}\subset A^{\Z}$ and consider the image of this cylinder set under the quotient map to $A^{\Z}/G$: $\overline{[aa'] _{(0,1)}} = \{(h b_i)_{i \in \Z} : b_0 = a, b_1 = a', h \in G\}\subset A^{\Z}/G$.  We have $\bar{\mu}(\overline{[aa']_{(0,1)}})  > 0$. Indeed, let $G' = \mathrm{Stab}(a) \cap \mathrm{Stab}(a')$. Then $G' \curvearrowright [aa']_{(0,1)}$ and $G'$ is finite by \pcom{the} properness of the action of~$G$ on~$\Delta$.  Since $G'$ is finite, we have a measurable fundamental domain $F$ for the action of $G'$ on $[aa']_{(0,1)}$ (\cite[Exercise 7.1.1]{Gao} and \cite[Exercise 7.1.6]{Gao}). We have that $\mu(F) = \frac{1}{\vert G' \vert} \mu([aa']_{(0,1)}) > 0$. Thus, $\bar{\mu}(\overline{[aa']_{(0,1)}})  = \mu(F) > 0$.

		 Let $Y = \overline{[aa']_{(0,1)}} \setminus \{(hb_i)_{i \in \Z} : h \in G \text{ and } p(b_j\to b_{j+1}) = 0 \text{ for some }j \in \Z \}$. Since $A$ is countable (since $\Delta$ is countable), we have that $\bar{\mu}(Y) = \bar{\mu}(\overline{[aa']_{(0,1)} })$.  Note that  all elements of $Y$ are then of the form $(h b_i)_{i \in \Z}$, where for each $j$, we have that $p(b_j \to b_{j+1}) > 0$, so that $b_j = (c, \ell)$ and $b_{j+1} = (c', \ell + 1)$, and $c,c'$ share a panel of type $s_{\ell}$.

		\textcolor{black}{Recalling from Lemma \ref{3.51} that $\bar{\mu}$ is a finite measure}, by Poincaré recurrence (see, e.g.\ \cite[Theorem~1.4]{Walters}) applied to the set $Y$ and the shift map $T \curvearrowright A^{\Z}/G$, we have that there exist $n > 0$ and some $(h b_i)_{i \in \Z} \in Y$ such that $T^{n} ((h b_i)_{i \in \Z}) \in Y$. Lifting back up to $A^{\Z}$, we obtain a sequence $(a_i)_{i \in \Z} \in [aa']_{(0,1)} $ such that $p(a_j \to a_{j+1}) > 0$ for all $j$ and such that $a_{n} = ga$ for some $g \in G$. 
		\end{proof}

	Note that in 
 Proposition \ref{3.6}, we have that $n = 0$ modulo $2k$ since $p(a_j \to a_{j+1}) > 0$ for all $j$.
	
	\vs

	\textit{Conclusion of the proof of Proposition \ref{3.11}}:
	
	\vs 
	
	Let $c_1$ be the first chamber of the \textcolor{black}{above geodesic gallery} $\gamma$ between $\Omega$ and $\Omega'$ and let $\sigma$ be the panel of type $s_0$ in $\Omega$ containing $c_1$. By Lemma \ref{3.5}, since $\Omega$ has a branching panel, every panel in $\Omega$ branches, so $\sigma$ branches.  Let  $c_2, c_3$ be two other distinct chambers in $\sigma$.
	
\begin{figure}[H]
	\centering
	\includegraphics[width=0.54\linewidth]{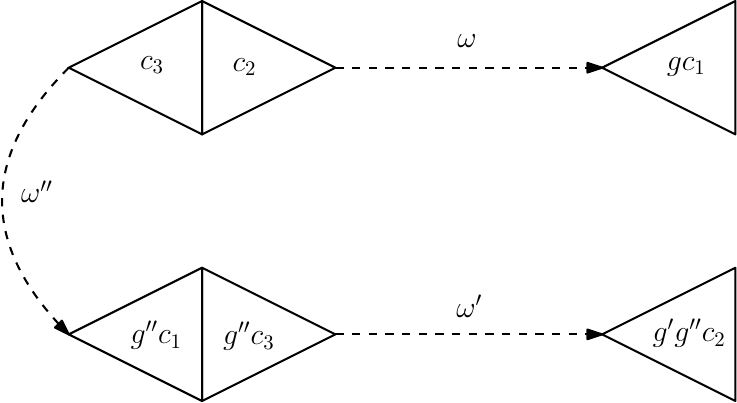}
	\caption{The \pcom{gallery $\omega^{-1}\omega''\omega'$ } produced from the branching panel $\sigma$.}
	\label{fig:dumbbell-graph}
\end{figure}

	 Apply Proposition \ref{3.6} to produce a sequence of states $((d_1,1),\ldots, (d_n, 0))$ whose chambers $d_i$ form a gallery $\omega$ from $c_2$ to  $gc_1$ for some $g \in G$. Similarly, produce a sequence of states whose chambers form a gallery $\omega''$ from $c_3$ to $g''c_1$ for some $g'' \in G$ and a sequence of states whose chambers form a gallery $\omega'$ from $g''c_3$ to $g'g''c_2$ for some $g' \in G$; see Figure~\ref{fig:dumbbell-graph}. Note that in the sequences produced by Proposition \ref{3.6}, adjacent states have positive transition probability, hence each of $\omega, \omega', \omega''$ has type of the form $s_0(s_0\cdots s_{2k-1})^m$ for some $m \in \N$, and hence is a geodesic gallery by Lemma \ref{3.2}. 

	 
 \pcom{The next claim establishes that the free group on $g,g'$ acts freely on $\langle g,g'\rangle \sigma$. (In fact, $\langle g,g'\rangle \omega^{-1}\omega''\omega'$ can be thought of as a ``tree of chambers'', but we will not need this more general statement.)}

\textit{Claim:}  Let $\sigma$ be the initial branching panel above in the wall $\Omega$. Then for any non-trivial freely reduced word $u$ in $g,g'$, we have $u\sigma \neq \sigma$ in $\Delta$. 
		
\begin{figure}[H]
	\centering
	\includegraphics[width=0.94\linewidth]{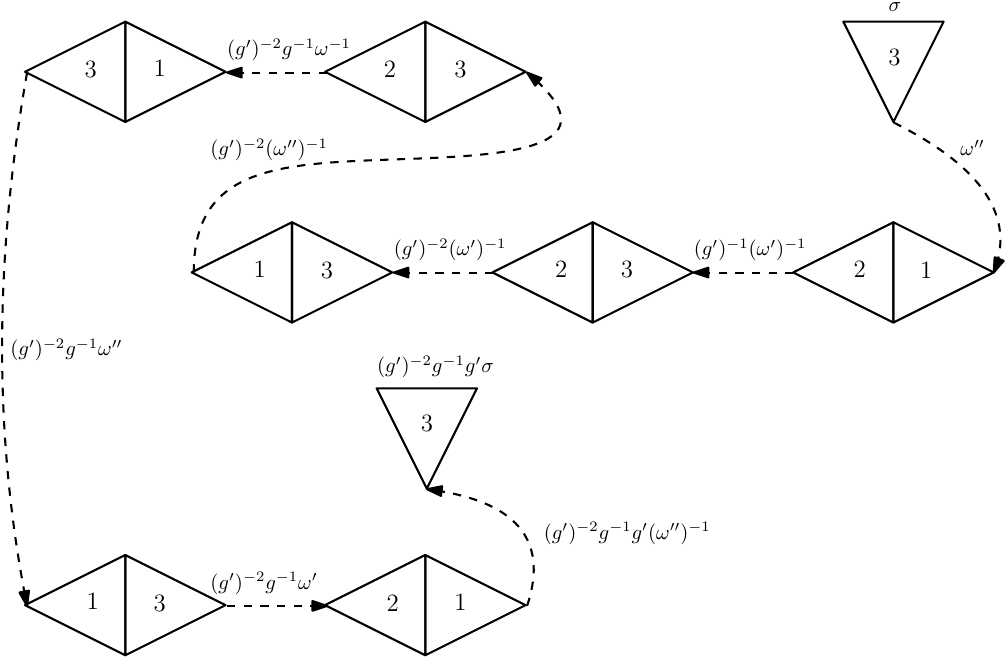}
	\caption{An example of a gallery joining $\sigma$ and $(g')^{-2}g^{-1} g'\sigma$. The numbers on the chambers indicate of which $c_i$ they are translates.}
	\label{fig:example-3}
\end{figure}

\begin{proof}[Proof of claim]
	
	Write $u$ as a word $u = u_1 u_2 \cdots u_m$, where the $u_i$ are alternating powers of $g$ and $g'$. Denote $u(i) = u_1 u_2 \cdots u_i$ for each $1 \leq i \leq m$ and let $u(0) = 1$. We show that for each $i$, we can connect  $u(i-1) \sigma$ to $u(i) \sigma$ with a gallery $\omega_i$ satisfying the following: 
	
	\begin{enumerate}[label = (\alph*)]
		\item $\omega_i$ is a concatenation of $\langle g,g' \rangle$-translates of $\omega$, $\omega'$ and $\omega''$,
		\item $\omega_i$ has type of the form $s_0w^{n_i}$ for some $n_i \in \N$, (recall that $w = s_0 \cdots s_{2k-1}$),
		\item the ending chamber of $\omega_i$ is different from the starting chamber of $\omega_{i+1}$.
	\end{enumerate}
	
	We first show that each $\omega_i$ satisfies (a) and (b). We consider the following cases:
	
	\begin{enumerate}[label = (\roman*)]
		\item $u_{i} = g^n$ for some $n \in \Z \setminus \{0\}$. Then $u(i) \sigma = u(i-1) g^n\sigma$. We can connect $\sigma$ to $g^n \sigma$ by $\bigcup_{j=0}^{n-1} g^j \omega$ if $n > 0$ or $\bigcup_{j=1}^{-n} g^{-j} \omega^{-1}$ if $n < 0$. Therefore, we set $\omega_i = u(i-1) \bigcup_{j=0}^{n-1} g^j \omega$ if $n > 0$ and $\omega_i = u(i-1) \bigcup_{j=1}^{-n} g^{-j} \omega^{-1}$ if $n < 0$. Thus, in this case we have that the type of $\omega_i$ is of the form $s_0w^{n_i}$, since the type of $\omega$ is of this form and the starting and ending chamber of~$\omega$ are different. See Figure~\ref{fig:example-1} for an illustration of an example. Note that the type of $\sigma$ (and hence all of its translates) is $s_0$.

\begin{figure}[H]
	\centering
	\includegraphics[width=0.9\linewidth]{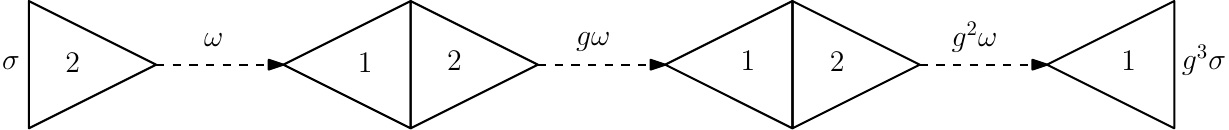}
	\caption{An example of a gallery joining $\sigma$ and $g^3 \sigma$. The numbers on the chambers indicate of which $c_i$ they are translates.}
	\label{fig:example-1}
\end{figure}
		
		\item $u_i = (g')^{n}$ for some $n \in \Z \setminus \{0\}$.  Then $u(i) \sigma = u(i-1) (g')^n\sigma$. For $n > 0$, we can connect $\sigma$ to $(g')^n \sigma$ by the concatenation  $\omega'' \cup (\bigcup_{j=0}^{n-1}(g')^j \omega') \cup (g')^n(\omega'')^{-1}$ and if $n < 0$, we can connect $\sigma$ to $(g')^n \sigma$ by the concatenation $\omega'' \cup (\bigcup_{j=1}^{-n} (g')^{-j} (\omega')^{-1}) \cup (g')^n(\omega'')^{-1} $, which has type of the form $s_0w^{n_i}$ for some $n_i \in \N$ since each of $\omega'', \omega'$ has type of this form and the starting and ending chamber of $\omega'$ and the ending chamber of $\omega''$ are distinct. Thus, $u(i-1) \sigma$ and $u(i) \sigma$ are connected by
        $$\omega_i=u(i-1) (\omega'' \cup (\bigcup_{j=0}^{n-1}(g')^j \omega') \cup (g')^n(\omega'')^{-1}) \; \mathrm{if} \; n>0,\; \mathrm{or}
        $$
        $$\omega_i=u(i-1) (\omega'' \cup (\bigcup_{j=1}^{-n} (g')^{-j} (\omega')^{-1}) \cup (g')^n(\omega'')^{-1}) \; \mathrm{if} \; n<0,
        $$
  which therefore have labels of the desired form $s_0w^{n_i}$ for some $n_i \in \N$. See Figure~\ref{fig:example-2} for an illustration of an example. 
		
\begin{figure}[H]
	\centering
	\includegraphics[width=0.75\linewidth]{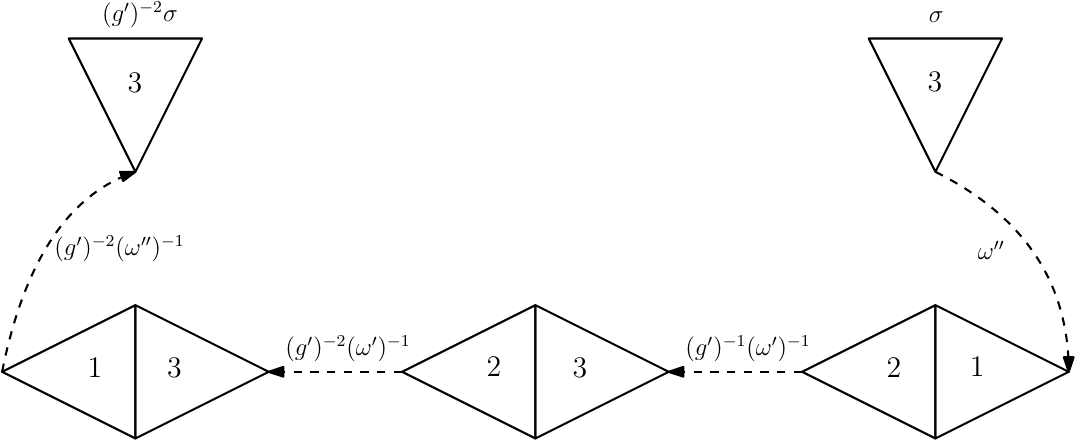}
	\caption{An example of a gallery joining $\sigma$ and $(g')^{-2} \sigma$. The numbers on the chambers indicate of which $c_i$ they are translates.}
	\label{fig:example-2}
\end{figure}
	\end{enumerate}

	Now we show that the ending chamber of $\omega_i$ is different from the starting chamber of $\omega_{i+1}$. 
	
	By the cases (i) and (ii) above, either for some $h \in G$, the ending chamber of $\omega_i$ is of the form $hc_1$ or $hc_2$ and the starting chamber of $\omega_{i+1}$ is of the form $hc_3$ (when $u_i$ is a power of $g$ and $u_{i+1}$ is a power of $g'$),  or for some $h \in G$, the ending chamber of $\omega_i$ is of the form $hc_3$ and the starting chamber of $\omega_{i+1}$ is of the form $hc_1$ or $hc_2$ (when $u_i$ is a power of $g'$ and $u_{i+1}$ is a power of $g$). 	Therefore, the $\omega_i$ satisfy (c). 
	
	Thus, the type of each $\omega_i$ is of the form $s_0 w^n$ and the starting and ending chamber of $\omega_i$ and~$\omega_{i+1}$ are distinct. Therefore, letting $\gamma = \bigcup_{i=1}^m \omega_i$ be the concatenation of the $\omega_i$ galleries, we have that $\gamma$ has type of the form $s_0w^{n_1 + n_2 + \cdots  + n_m}$. By Lemma \ref{3.2}, $s_0w^{n_1 + n_2 + \cdots  + n_m}$ is a geodesic word, and so we have that $\gamma$ is a geodesic gallery in $\Delta$. Therefore, $\gamma$ has distinct endpoints, and so $\sigma\neq u\sigma$.
    In Figure~\ref{fig:example-3}, see an example of the gallery $\gamma$ for $u = (g')^{-2}g^{-1} g'$. 
    \end{proof} 

By the above claim, we obtain that $\langle g, g' \rangle \cong F_2$. Therefore, we have that $G$ contains  $F_2$ as a subgroup, concluding the proof of Proposition \ref{3.11}. \hfill $\qed$
  
\section{Appendix}

\pcom{We outline the proof of Theorem \ref{B}} illustrating how the weak Tits alternative of a group can be deduced from its action on a product of CAT(0) spaces satisfying rank rigidity. This theorem and its proof were communicated to us by Pierre-Emmanuel Caprace.

First, we recall the notion of rank rigidity and irreducibility of metric spaces.

A metric space $X$ is \textbf{irreducible} if in each decomposition $X = Y \times Z$, one of the factors $Y$ or~$Z$ is a point, \textcolor{black}{where we define the product metric $d$ on $X$ in terms of the metrics $d_Y, d_Z$ on $Y,Z$ via $d((y_1,z_1), (y_2,z_2)) = \sqrt{d_Y(y_1,y_2)^2 + d_Z(z_1,z_2)^2}$. }

A \textbf{flat} in a metric space $X$ is an isometrically embedded copy of $\R^n$ in $X$. A \textbf{half-flat} in $X$ is an isometrically embedded half of $\R^n$ in $X$.

Given a CAT(0) space $X$, a hyperbolic element $g \in \mathrm{Isom}(X)$ is a \textbf{rank one isometry} if $g$ does not have a geodesic axis in $X$ which bounds a half-flat in $X$. 

\begin{defn}
    An irreducible \textcolor{black}{$\mathrm{CAT(0)}$} space $X$ satisfies \textbf{rank rigidity} if $X$ is either a Euclidean building, a symmetric space of non-compact type or $\mathrm{Isom}(X)$ contains a rank one isometry. 
\end{defn}


\begin{proof}[\pcom{Proof of Theorem \ref{B}}]

    First, if some $X_i$ is quasi-isometric to $\R$, then we can replace $X_i$ by $\R$ in the product decomposition of $X$. Then grouping together all of the $\R$ factors in $X$, we can rewrite $X$ in the form $X = X_1 \times \cdots \times X_k \times \R^m$ for $k \leq n$ and $m \geq 0$, where each $X_i$ is an irreducible metric space that is neither a point nor quasi-isometric to $\R$. 

    By \cite[Corollary 1.3]{MR2399098}, up to passing to a finite index subgroup of $G$,  we can assume that the image of $G$ in $\mathrm{Isom}(X)$ lies in $\mathrm{Isom}(X_1) \times \cdots \times \mathrm{Isom}(X_k) \times \mathrm{Isom}(\R^m)$. For each $i = 1,\ldots,k$, let $G_i$ denote the image of $G$ in $\mathrm{Isom}(X_i)$. 

    We split into cases based on the type of each $X_i$ from rank rigidity. 

    \begin{enumerate}[label = (\roman*)]
			\item Each $X_i$ is a Euclidean building. Then $X$ is of the form $X = B \times \R^m$, where $B$ is a \pcom{(possibly reducible)} Euclidean building. \pcom{Arguing as in the proof of \cite[Theorem F]{BallmanBrin}, one can prove that $G$ satisfies the weak Tits alternative.}
			\item Some $X_i$ is a symmetric space of non-compact type. We have that $G_i$ is linear, by \cite[Proposition 2.1.1]{MR1441541}. Since symmetric spaces are geodesically complete (i.e.\ every geodesic segment extends to a geodesic line), by \cite[Lemma II.6.20]{BH99}, we have that the action of $G_i$ on $X_i$ is minimal (i.e.\ $X_i$ does not contain any closed, convex $G_i$-invariant subspaces other than $\emptyset$ and $X_i$). Hence, by \cite[Theorem L]{MR3072802}, $G_i$ does not fix any point in the boundary $\partial X_i$ of $X_i$. Also, there is no $G_i$-invariant flat in $X_i$, since otherwise by \pcom{the} minimality of the action of $G_i$ on $X_i$, such a flat would equal $X_i$, contradicting the fact that $X_i$ is irreducible and not equal to $\R$.
            
            By the Adams--Ballmann Theorem \cite{MR1645958} applied to $G_i < \mathrm{Isom}(X_i)$, we then have that $G_i$ is non-amenable. Therefore by the classical Tits alternative for linear groups, we have that $F_2 < G_i$ and hence $F_2 < G$.  
            \item Some $X_i$ has a rank one isometry. Since $G_i$ acts cocompactly on $X_i$, we have that the limit set of $G_i$ in $\partial X_i$ is all of $\partial X_i$. By \cite[Proposition 3.5]{MR3449592}, we have that $G_i$ contains a rank one isometry. 
            
            We will show that $G_i$ does not fix a point in the boundary $\partial X_i$ of $X_i$. Suppose that $G_i$ did fix a point in $\partial X_i$. Then by \cite[Proposition 3.15]{MR2574741}, there exists a closed, convex subspace $Y \subseteq X_i$ with $\partial Y = \partial X_i$ and such that $Y$ splits $G_i$-equivariantly as $Y \cong \R^n \times W$, and any $G_i$-fixed point in the boundary lies in $\partial \R^n$. Since we assumed that $G_i$ has a fixed point in $\partial X_i$, we have $n > 0$. 
            
            We will show that \textcolor{black}{every} axis $\gamma$ of any rank one isometry in $G_i$ lies in a copy of $\R \times \R$ or $\R \times [0, \infty)$ in $Y$, hence $\gamma$ bounds a half-flat in $X_i$, contradicting the definition of rank one isometry. We consider the following cases: 
            
            \begin{itemize}
                \item Suppose that $\gamma$ does not lie in any factor of $Y$. Then the projection of $\gamma$ onto each factor of $Y$ is a geodesic line, and so $\gamma \subset \R \times \R$. 
                \item Suppose that $\gamma$ lies in the $\R^n$ factor of $Y$. Then we have $n = 1$, since otherwise $\gamma$ clearly bounds a half-flat. Since $W$ is unbounded (as otherwise $Y$, and hence $X_i$, would be quasi-isometric to $\R$, contradicting our assumption at the beginning of the proof), we have that $W$ contains a geodesic ray. Hence, $\gamma \subset \gamma \times [0, \infty) \cong \R \times [0, \infty)$. 
                \item Suppose that $\gamma$ lies in $W$. Then since $n \geq 1$, we have $\gamma \subset \R \times \gamma \cong \R \times \R$. 
            \end{itemize}
            
            Thus, we obtain a contradiction in each case, and so we conclude that $G_i$ does not fix a point in $\partial X_i$.
            
            Since $G_i$ does not fix any point of $\partial X_i$, by \cite[Proposition 3.4]{MR2585575} we have that either $G_i$ contains two independent rank one isometries or $G_i$ stabilizes some geodesic line in $X_i$. The latter option cannot occur since then by \pcom{the} cocompactness of the action of $G_i$ on $X_i$, we would have that $X_i$ is quasi-isometric to $\R$, contradicting our assumption made at the beginning of the proof. 
            
            Therefore, we have that $G_i$ contains two independent rank one isometries, whose powers therefore generate a copy of $F_2$ in $G_i$ by \cite[Proposition 3.4(2)]{MR2585575}. Hence, $F_2 < G$. 
		\end{enumerate}

\end{proof}

As an immediate consequence of Theorem \ref{B}, we obtain an alternate proof of \textcolor{black}{Theorem \ref{A}}.

\begin{proof}[Proof of \textcolor{black}{Theorem \ref{A}} using \pcom{Theorem \ref{B}}]
    
    Let $X$ be a building with a group $G$ acting geometrically on $X$ by type-preserving automorphisms.   Decompose $X$ as $X \cong X_1 \times \cdots \times X_n$ for $X_i$ irreducible buildings (this decomposition corresponds to the decomposition of the underlying Coxeter group~$W$ into irreducible factors: $W \cong W_1 \times \cdots \times W_n$). By Lemma \ref{3.9}, we can assume that no $W_i$ is finite. \textcolor{black}{By \cite[Proposition 3.2]{MR3276856} each $X_i$ is either a Euclidean building or $\mathrm{Isom}(X_i)$ contains a rank one isometry.} By \pcom{Theorem \ref{B}}, we conclude that $G$ satisfies the weak Tits alternative. 
\end{proof} 

	\newpage
	\bibliographystyle{plain} 
	\bibliography{refs2}

\begin{thebibliography}{10}

\bibitem{BrownAbr}
Peter Abramenko and Kenneth~S. Brown.
\newblock {\em Buildings}, volume 248 of {\em Graduate Texts in Mathematics}.
\newblock Springer, New York, 2008.
\newblock Theory and applications.

\bibitem{MR1645958}
Scott Adams and Werner Ballmann.
\newblock Amenable isometry groups of {H}adamard spaces.
\newblock {\em Math. Ann.}, 312(1):183--195, 1998.

\bibitem{BallmanBrin}
Werner Ballmann and Michael Brin.
\newblock Orbihedra of nonpositive curvature.
\newblock {\em Inst. Hautes \'{E}tudes Sci. Publ. Math.}, (82):169--209 (1996), 1995.

\bibitem{MR0693651}
Hyman Bass and Alexander Lubotzky.
\newblock Automorphisms of groups and of schemes of finite type.
\newblock {\em Israel J. Math.}, 44(1):1--22, 1983.

\bibitem{BH99}
Martin~R. Bridson and Andr\'{e} Haefliger.
\newblock {\em Metric spaces of non-positive curvature}, volume 319 of {\em Grundlehren der mathematischen Wissenschaften [Fundamental Principles of Mathematical Sciences]}.
\newblock Springer-Verlag, Berlin, 1999.

\bibitem{BrinkB}
Brigitte Brink and Robert~B. Howlett.
\newblock A finiteness property and an automatic structure for {C}oxeter groups.
\newblock {\em Math. Ann.}, 296(1):179--190, 1993.

\bibitem{MR2110465}
Pierre-Emmanuel Caprace.
\newblock The thick frame of a weak twin building.
\newblock {\em Adv. Geom.}, 5(1):119--136, 2005.

\bibitem{MR2585575}
Pierre-Emmanuel Caprace and Koji Fujiwara.
\newblock Rank-one isometries of buildings and quasi-morphisms of {K}ac-{M}oody groups.
\newblock {\em Geom. Funct. Anal.}, 19(5):1296--1319, 2010.

\bibitem{MR3449592}
Pierre-Emmanuel Caprace and David Hume.
\newblock Orthogonal forms of {K}ac-{M}oody groups are acylindrically hyperbolic.
\newblock {\em Ann. Inst. Fourier (Grenoble)}, 65(6):2613--2640, 2015.

\bibitem{MR2574741}
Pierre-Emmanuel Caprace and Nicolas Monod.
\newblock Isometry groups of non-positively curved spaces: discrete subgroups.
\newblock {\em J. Topol.}, 2(4):701--746, 2009.

\bibitem{MR3072802}
Pierre-Emmanuel Caprace and Nicolas Monod.
\newblock Fixed points and amenability in non-positive curvature.
\newblock {\em Math. Ann.}, 356(4):1303--1337, 2013.

\bibitem{TwistRig}
Pierre-Emmanuel Caprace and Piotr Przytycki.
\newblock Twist-rigid {C}oxeter groups.
\newblock {\em Geom. Topol.}, 14(4):2243--2275, 2010.

\bibitem{BipolarCox}
Pierre-Emmanuel Caprace and Piotr Przytycki.
\newblock Bipolar {C}oxeter groups.
\newblock {\em J. Algebra}, 338:35--55, 2011.

\bibitem{MR3276856}
Corina Ciobotaru.
\newblock The flat closing problem for buildings.
\newblock {\em Algebr. Geom. Topol.}, 14(5):3089--3096, 2014.

\bibitem{Davis}
Michael~W. Davis.
\newblock {\em The geometry and topology of {C}oxeter groups}, volume~32 of {\em London Mathematical Society Monographs Series}.
\newblock Princeton University Press, Princeton, NJ, 2008.

\bibitem{DavisShapiro}
Michael~W. Davis and Michael~D. Shapiro.
\newblock Coxeter groups are automatic.
\newblock 1991.
\newblock Preprint.

\bibitem{MR1441541}
Patrick~B. Eberlein.
\newblock {\em Geometry of nonpositively curved manifolds}.
\newblock Chicago Lectures in Mathematics. University of Chicago Press, Chicago, IL, 1996.

\bibitem{MR2399098}
Thomas Foertsch and Alexander Lytchak.
\newblock The de {R}ham decomposition theorem for metric spaces.
\newblock {\em Geom. Funct. Anal.}, 18(1):120--143, 2008.

\bibitem{Gao}
Su~Gao.
\newblock {\em Invariant descriptive set theory}, volume 293 of {\em Pure and Applied Mathematics (Boca Raton)}.
\newblock CRC Press, Boca Raton, FL, 2009.

\bibitem{BuildingsApp}
Lizhen Ji.
\newblock Buildings and their applications in geometry and topology.
\newblock {\em Asian J. Math.}, 10(1):11--80, 2006.

\bibitem{Noskov2002}
Guennadi~A. Noskov and \`Ernest~B. Vinberg.
\newblock Strong {T}its alternative for subgroups of {C}oxeter groups.
\newblock {\em J. Lie Theory}, 12(1):259--264, 2002.

\bibitem{MR4305241}
Damian Osajda and Piotr Przytycki.
\newblock Tits alternative for groups acting properly on 2-dimensional recurrent complexes.
\newblock {\em Adv. Math.}, 391:Paper No. 107976, 22, 2021.
\newblock With an appendix by J. McCammond, Osajda and Przytycki.

\bibitem{osajda_przytycki_2022}
Damian Osajda and Piotr Przytycki.
\newblock Tits alternative for 2-dimensional {$\rm CAT(0)$} complexes.
\newblock {\em Forum Math. Pi}, 10:Paper No. e25, 19, 2022.

\bibitem{Ronan}
Mark Ronan.
\newblock {\em Lectures on buildings}.
\newblock University of Chicago Press, Chicago, IL, 2009.
\newblock Updated and revised.

\bibitem{SageevWise}
Michah Sageev and Daniel~T. Wise.
\newblock The {T}its alternative for {${\rm CAT}(0)$} cubical complexes.
\newblock {\em Bull. London Math. Soc.}, 37(5):706--710, 2005.

\bibitem{Tits1}
Jacques Tits.
\newblock Free subgroups in linear groups.
\newblock {\em J. Algebra}, 20:250--270, 1972.

\bibitem{Walters}
Peter Walters.
\newblock {\em An introduction to ergodic theory}, volume~79 of {\em Graduate Texts in Mathematics}.
\newblock Springer-Verlag, New York-Berlin, 1982.

\bibitem{MR1783960}
Burkhard Wilking.
\newblock On fundamental groups of manifolds of nonnegative curvature.
\newblock {\em Differential Geom. Appl.}, 13(2):129--165, 2000.

\end{thebibliography}

\end{document}